\documentclass[12pt]{article}
\usepackage[top=2.54cm, bottom=2.54cm, left=3.18cm, right=3.18cm]{geometry}
\usepackage{amsmath}
\usepackage{amssymb}
\usepackage{amsthm}
\usepackage{fancyhdr}
\usepackage{latexsym}
\usepackage{mathrsfs}
\usepackage{wasysym}
\usepackage{fancyhdr}
\usepackage{float}
\usepackage{graphicx}
\usepackage[compress]{cite}
\usepackage{subfigure}
\usepackage{comment}
\allowdisplaybreaks[4]

\newtheorem{theorem}{\bf Theorem}[section]

\newtheorem{lemma}[theorem]{Lemma}

\newtheorem{conjecture}[theorem]{Conjecture}

\begin{document}

\title{A bound on judicious bipartitions of directed graphs}{}

\author{Jianfeng Hou$^1$\thanks{Partially supported by NSFC grant
    11671087},~  Huawen Ma$^1$, ~ Xingxing Yu$^2$\thanks{Partially
    supported by NSF grant DMS-1600738 and the Hundred  Talents Program of Fujian Province}, ~Xia Zhang$^3$\thanks{Partially supported by the Shandong Provincial Natural Science Foundation, China (Grant No. ZR2014JL001) and the Excellent Young Scholars Research
Fund of Shandong Normal University of China.
} \\
{\small $^1$Center for Discrete Mathematics, Fuzhou University, Fujian, China, 350116}\\
{\small $^2$School of Mathematics, Georgia Institute of Technology, Atlanta, GA 30332-0160, USA}\\
{\small $^3$School of Mathematics and Statistics, Shandong Normal University, Shandong, China, 250014}
}
\date{}

\maketitle

\begin{abstract}
Judicious partitioning problems on graphs ask for partitions
that bound several quantities simultaneously, which have received a lot of
attentions lately.
Scott asked the following natural question: What is the maximum constant $c_d$ such that every
directed graph $D$ with $m$ arcs and minimum outdegree $d$ admits a bipartition $V(D)= V_1\cup V_2$
satisfying $\min\{e(V_1, V_2), e(V_2, V_1)\}\ge c_d m$? Here, for
$i=1,2$, $e(V_{i},V_{3-i})$ denotes the number of arcs in $D$ from
$V_{i}$ to $V_{3-i}$.   Lee,  Loh, and Sudakov
conjectured that  every directed graph $D$ with $m$ arcs and minimum
outdegree at least $d\ge 2$ admits a bipartition $V(D)=V_1\cup V_2$
such that
\[
\min\{e(V_1,V_2),e(V_2,V_1)\}\geq \Big(\frac{d-1}{2(2d-1)}+ o(1)\Big)m.
\]
We show that this conjecture holds under the
additional natural condition that the minimum indegree is also at least $d$.
\end{abstract}



\section{Introduction}

The Max-Cut problem asks for a cut with maximum weight in a weighted
graph, which has been studied extensively from both
combinatorial and computational perspectives. It is  NP-hard
even when restricted to the class of triangle-free cubic graphs
\cite{Yannakakis1978}.  A simple calculation shows that every graph with
$m$ edges has a cut with at least $m/2$ edges. Answering a question of Erd\H{o}s, Edwards
\cite{Edwards1973,Edwards1975} improved this lower bound to
$m/2+(\sqrt{2m+1/4}-1/2)/4$, which is tight for complete
graphs with odd order. Moreover, for certain range of $m$, Alon \cite{Alon1996} gave an additive improvement of order $m^{1/4}$. For special classes of graphs, such as subcubic graphs \cite{Xu2008, Zhu2009}, the main term in the Edwards' bound may be improved.

\medskip

For a directed graph (digraph for short)  $D$ and $S,T\subseteq V(D)$, we use $e(S,T)$ to denote
the number of arcs of $D$ directed from $S$ to $T$.
A directed (unweighted) version of the Max-Cut problem is to find a partition
$V(D)=V_1\cup V_2$ of a given digraph $D$ that maximizes  $e(V_1,
V_2)$. It is easy to see from the above bound of Edwards that every
digraph $D$ with $m$ arcs has a partition $V(D)=V_1\cup V_2$ with
$e(V_1,V_2)\ge m/4+(\sqrt{2m+1/4}-1/2)/8$, and that any regular
orientation of a complete graph of odd order shows that this bound is
tight. Alon,  Bollob\'as, Gy\'arf\'as, Lehel, and  Scott
\cite{Alon2007} considered partitions of acyclic digraphs and proved
that every acyclic digraph $D$  with $m$ arcs has a  partition
$V(D)=V_1\cup V_2$ such that $e(V_1, V_2)\ge m/4 +
\Omega(m^{2/3})$. For positive integers $k$ and $l$, let
$\mathcal{D}(k,l)$ denote the family of digraphs such that every
vertex has indegree at most $k$ or outdegree at most $l$. It is proved in  \cite{Alon2007} that
every digraph $D\in \mathcal{D}(k,l)$ with $m$ arcs has a  partition $V(D)=V_1\cup V_2$ such that $e(V_1, V_2)\ge (k+l+2)m/(4k+4l+6)$. For digraphs $D\in \mathcal{D}(1,1)$ with $m$ arcs, Chen, Gu, and Li \cite{Chen2014} showed that $D$ has a  partition $V(D)=V_1\cup V_2$ such that $e(V_1, V_2)\ge 3(m-1)/8$.

\medskip

In practice, one often needs to find a partition of a given graph or
digraph that  simultaneously bounds several quantities. Such problems are called $Judicious$ $Partitioning$ $Problems$ by Bollob\'{a}s and Scott \cite{Bollobas2002_1}. In \cite{Scott2005}, Scott ask the following natural question: What is the maximum constant $c_d$ such that every
digraph $D$ with $m$ arcs and minimum outdegree $d$ admits a bipartition $V(D)= V_1\cup V_2$
satisfying
\[
\min\{e(V_1, V_2), e(V_2, V_1)\}\ge c_d m?
\]

For $d=1$, the graph, obtained from $K_{1,n-1}$ by adding a single edge inside
the part of size $n-1$, admits an orientation in which the minimum outdegree is 1 and $\min\{e(V_1, V_2), e(V_2, V_1)\}\le 1$ for every partition
$V(D)= V_1\cup V_2$. Hence,  $c_1=0$. Lee, Loh, and Sudakov
\cite{Lee2014} studied Scott's question and made the following
\begin{conjecture}[Lee, Loh, and Sudakov \cite{Lee2014}]\label{conjecture} \quad
Let $d$ be an integer satisfying $d\geq 2$. Every digraph $D$ with $m$ arcs and minimum outdegree at least $d$ admits a bipartition $V(D)=V_1\cup V_2$ with
$$
\min\{e(V_1,V_2),e(V_2,V_1)\}\geq \Big(\frac{d-1}{2(2d-1)}+ o(1)\Big)m.$$
\end{conjecture}
In the same paper, Lee, Loh, and Sudakov verify Conjecture~\ref{conjecture} for $d=2,3$, and  noted that
their method is not adequate for  $d\geq 4$.

\medskip
It is not clear to us whether or not the minimum outdegree condition
alone is sufficient for Conjecture~\ref{conjecture} with large $d$. However, we show in this paper  that
Conjecture \ref{conjecture} holds under the natural (additional) assumption that
the minimum indegree of $D$ is at least $d$.

\begin{theorem}\label{main1} \quad
Let $D$ be a digraph with $m$ arcs and assume that both the
minimum outdegree  and the minimum indegree of $D$ are at least
$d\ge 2$. Then $D$ admits a bipartition $V(D)=V_{1}\cup V_{2}$ with
$$\min\{e(V_{1},V_{2}),e(V_{2},V_{1})\}\geq\Big(\frac{d-1}{2(2d-1)}+o(1)\Big)m.$$
\end{theorem}

In the remainder of this section, we describe notation and
terminology used in the proof of Theorem~\ref{main1}. In Section 2, we
collect previous results (as well as a concentration inequality) that we will need. The proof of
Theorem~\ref{main1} will be given in Section 3. In Section 4, we offer
some concluding remarks.

\medskip

All digraphs  considered in this paper are finite with no loops and no parallel arcs. For a digraph $D$,
we denote by $V ( D )$ and $A(D)$ the $vertex$ $set$ and the $arc$ $set$ of $D$, respectively.
Let $|D|:=|V(D)|$ be the \emph{order} of $D$, and $e(D):=|A(D)|$ be the
$size$ of $D$. The \emph{underlying graph} of $D$ is obtained from $D$
by ignoring arc orientations and removing
redundant parallel edges.

Given two vertices $x,y$ in a digraph $D$, we write $xy$ for the arc directed from $x$ to $y$.
Let $N^+_D(x):=\{z: xz\in A(D)\}$ and $N^-_D ( x ):=\{z: zx\in
A(D)\}$.
Then $d^+_D( x ):= | N^+_D ( x )|$ and $d^-_D( x ):= | N^-_D ( x )|$ are
the the outdegree and indegree of $x$, respectively.
The minimum outdegree of $D$ is  $\delta^+(D):=\min\{d^+_D(x): x\in V(D)\}$ and the minimum  indegree of $D$ is  $\delta^-(D):=\min\{d^-_D(x): x\in V(D)\}$.
Let $\delta^0(D):=\min\{\delta^+(D),\delta^-(D)\}$.
The \emph{degree } of $x\in V(D)$ is  defined as $d_D(x) :=
d^+_D(x) + d^-_D(x)$. We use  $\delta(D) :=\min\{d_D(x) : x \in V (D)\}$
and $\Delta(D) := \max\{d_D(x) : x \in V (D)\}$ to denote
the \emph{minimum degree} and \emph{maximum degree} of $D$, respectively.

Let $D$ be a digraph. For $X\subseteq V(D)$,
the subgraph of $D$ \emph{induced} by $X$ is denoted by $D [X]$. Let $e_D(X)$ denote the number of arcs in $D [X]$, and
let $D-X$ denote the digraph obtained from $D$ by deleting $X$ and all arcs incident with $X$.

Throughout this paper, we drop the reference to $D$ in the above
notations if there is no danger of confusion.

\medskip

\section{Lemmas}
In this section, we list previous results needed in our
proof of Theorem \ref{main1}.
Lee, Loh, and Sudakov \cite{Lee2013} introduced the notion  of a tight component in a graph when analyzing
bisections of graphs. A connected graph $T$ is $tight$ if
\begin{itemize}
\item for every vertex $v\in V(T)$, $T- v$ contains a perfect matching, and

\item  for every vertex $v\in V(T)$ and every perfect matching $M$ of $T- v$, no edge in $M$ has exactly one end adjacent to $v$.
\end{itemize}
Note that if $T$ is tight, then the order of $T$ is odd and  each
vertex in $T$ has even degree.  Lu, Wang, and  Yu \cite{lu2015} gave
the following characterization of tight graphs.

\begin{lemma} [Lu, Wang, and  Yu \cite{lu2015}]\label{tight-graph-clique}
A connected graph $G$ is tight iff every block of $G$ is an odd clique.
\end{lemma}

The following result bounds the number of tight components in a graph,
which is an easy consequence of Lemma~\ref{tight-graph-clique}.

\begin{lemma}[Lee, Loh, and Sudakov \cite{Lee2013}]\label{Lee1}
For each integer $i\geq 0$, let $d_{i}$ denote the number of vertices with degree
$i$ in a graph $G$. Then the number of tight components $\tau$ of $G$ satisfies
\[
\tau\leq \frac{d_{0}}{1}+\frac{d_{2}}{3}+\frac{d_{4}}{5}+\cdots.
\]
\end{lemma}

We also need  a result from \cite{Hou2017}, which
provides a bound on $d_0/1+d_2/3+\cdots$ (and, hence, on $\tau$) under additional
constraints.

\begin{lemma}[Hou and Wu \cite{Hou2017}]\label{hou2}
Let $n$, $\rho$, $\alpha$, $\delta$ be fixed nonnegative integers with $\delta-\alpha\geq1$,
and let $d_{i}$ be a real number with $0\le d_i\le n$ for $i\in\{0,1,\cdots,n-1\}$ such that
$$\sum_{i=0}^{n-1} d_{i}\leq n \mbox{ and }
(\delta-\alpha)d_{0}+(\delta-\alpha-1)d_{1}+\cdots+d_{\delta-\alpha-1}\leq\rho.$$
Then
$$\frac{d_{0}}{1}+\frac{d_{2}}{3}+\cdots\leq\frac{n+\rho}{\delta-\alpha+1}.$$
\end{lemma}

We now turn to lemmas on partitions for digraphs, first of which
concerns dense digraphs.

\begin{lemma}[Lee, Loh, and Sudakov \cite{Lee2013}]\label{Lee2}
Let $D$ be a digraph with $n$ vertices and $m$ arcs. For any $\varepsilon>0$, if $m\geq8n/\varepsilon^{2}$
or $\Delta(D)\leq\varepsilon^{2}m/4$, then $D$ admits a bipartition $V(D)=V_{1}\cup V_{2}$ with
$\min\{e(V_{1},V_{2}),e(V_{2},V_{1})\}\geq m/4-\varepsilon m$.
\end{lemma}

One approach for finding a ``good'' bipartition of a digraph $D$ is to
first partition $V(D)$ into two sets $X$ and $Y$ (with $X$ consisting of certain high
degree vertices), then partition $X$ into two sets $X_1$ and $X_2$ with certain
property, and finally apply a randomized algorithm to partition $Y$.
The following lemma will be used to identify a useful partition of $X$.

\begin{lemma}[Hou, Wu, and Yan \cite{hou}]\label{hou1}
Let $D$ be a digraph and $V(D)=X\cup Y$ be a partition of $D$ with $e(X)=0$.
For any given partition $X=X_{1}\cup X_{2}$, define its gap to be
$$\theta(X_1,X_2)=\big(e(X_{1},Y)+e(Y,X_{2}))-(e(X_{2},Y)+e(Y,X_{1})\big).$$
If $X_{1}\cup X_{2}$ is a partition of $X$ which
minimizes $|\theta(X_1,X_2)|$, then $|\theta(X_1,X_2)|\leq|Y|$.
\end{lemma}
Note that $|\theta(X_1,X_2)|$ is the absolute value of
$\theta(X_1,X_2)$, while $|Y|$ is the cardinality of $Y$. We remark
that $\theta(X_1,X_2)$ was introduced by Lee, Loh and Sudakov \cite{Lee2014}.

\medskip

The next lemma allows us to extend a partial bipartition $X_1,X_2$ of a digraph to a ``good'' bipartition of the entire digraph. This will be used several
times in the proof of Theorem \ref{main1}.

\begin{lemma}[Lee, Loh, and Sudakov \cite{Lee2014}]\label{main-lemma}
For any real constants $C$, $\varepsilon>0$, there exist $\gamma$, $n_{0}>0$ for which the following holds.
Let $D$ be a digraph with $n\geq n_{0}$ vertices and at most $Cn$ arcs. Suppose $X\subseteq V(D)$ is a set
of at most $\gamma n$ vertices which have been partitioned into $X_{1}\cup X_{2}$.
Let $Y=V(D)\setminus X$ and $\tau$ be the number of tight components in $G[Y]$, where $G$ is the
underlying graph of $D$. If every vertex in $Y$ has degree at most
$\gamma n$ in $G$, then there is a bipartition
$V(D)=V_{1}\cup V_{2}$ with $X_{i}\subseteq V_{i}$ for $i=1,2$ such that
$$e(V_{1},V_{2})\geq e(X_{1},X_{2})+\frac{e(X_{1},Y)+e(Y,X_{2})}{2}+\frac{e(Y)}{4}+\frac{n-\tau}{8}-\varepsilon n,$$
$$e(V_{2},V_{1})\geq e(X_{2},X_{1})+\frac{e(X_{2},Y)+e(Y,X_{1})}{2}+\frac{e(Y)}{4}+\frac{n-\tau}{8}-\varepsilon n.$$
\end{lemma}

For the step to partition $Y$, we will need  the  following
version of the Azuma-Hoeffding inequality \cite{Azuma1967,Hoeffding1963}, see Corollary 2.27 in  Janson, {\L}uczak, and Ruci\'{n}ski \cite{janson2000}.

\begin{lemma}[Azuma \cite{Azuma1967}, Hoeffding \cite{Hoeffding1963}]\label{AH}
Let $Z_1,\ldots, Z_n$ be independent random variables taking values in $\{1,\ldots,k\}$, let $Z:= (Z_1,\ldots, Z_n)$, and let $f : \{1,\ldots,k\}^n$ $\rightarrow \mathbb{N}$ such that $|f (Y )- f (Y')|\leq c_i$ for any $Y , Y' \in\{1,\ldots,k\}^n$ which differ only in the $i$th coordinate. Then for any $z > 0$,
\[\mathbb{P}\big(f(Z)\geq \mathbb{E}\big(f(Z)\big)+z\big)\leq \exp\left(\frac{-z^2}{2\sum_{i=1}^nc_i^2}\right),\]
\[\mathbb{P}\big(f(Z)\leq \mathbb{E}\big(f(Z)\big)-z\big)\leq \exp\left(\frac{-z^2}{2\sum_{i=1}^nc_i^2}\right).\]
\end{lemma}

\section{Proof of Theorem \ref{main1}}

In this section, we prove  Theorem \ref{main1}. Our approach is
similar to that used in Lee, Loh, and Sudakov \cite{Lee2014}, which was
used first by Bollob\'as and Scott \cite{Scott1999}
and again by Ma, Yen, and Yu \cite{MaJCTB}.
However, additional ideas are needed to make this approach work.

Let $d\ge 2$ be an integer, let $D$ be an $n$-vertex digraph with $m$ arcs, and assume
$\delta^0(D):= d\ge 2$ (i.e., both minimum outdegree $\delta^+(D)$ and minimum
indegree $\delta^-(D)$ are at least $d$).
Then, since $\delta^+(D)\ge d$,  $$m=\sum_{v\in  V(D)}d^+(v)\ge dn.$$
We need to show that $D$
admits a bipartition $V(D)=V_{1}\cup V_{2}$ with
$$\min\{e(V_{1},V_{2}),e(V_{2},V_{1})\}\geq\Big(\frac{d-1}{2(2d-1)}+o(1)\Big)m,$$
and we will reduce this to proving \eqref{mian-ineq-pf}.
\medskip

We may  assume that $$d\ge 4,$$
as  Lee, Loh, and Sudakov \cite{Lee2014}  showed that, for $d\le 3$,
such a bipartition exists.
We may also assume that
\begin{equation}\label{m=on=delta}
m\le 128(2d-1)^{2}n \text{ and } \Delta(D)\ge \frac{m}{64(2d-1)^{2}},
\end{equation}
for, otherwise, applying Lemma \ref{Lee2} with $\varepsilon=1/(8d-4)$
gives a partition $V(D)=V_{1}\cup V_{2}$ such that
\[
\min\{e(V_{1},V_{2}),e(V_{2},V_{1})\}\geq\Big(\frac{1}{4}-\frac{1}{8d-4}\Big)m=\frac{d-1}{2(2d-1)}m.
\]
Thus, we may  assume that
$$\mbox{$n$ is sufficiently large.}$$

\medskip

For clarity of presentation, we divide the remainder of the proof into nine steps.

\medskip

{\bf Step 1. Identifying the set $X$,  and partitioning $X$ to $X_1\cup X_2$.}

Let $X$ be the set of vertices of $D$ with degree at least $n^{3/4}$, and let $Y=V(D)-X$ and $n':=|Y|$. Then
\[
2m\ge \sum_{v\in X}d(v)\ge |X|n^{3/4}.
\]
It follows from \eqref{m=on=delta} that
\begin{equation}\label{bound-X-n'}
|X|\le 256(2d-1)^{2}n^{1/4}=O(n^{1/4}).
\end{equation}
Thus, $e(X)\leq|X|^{2}=O(n^{1/2})$. Therefore, for the sake of
simplicity, we remove all the arcs of $D$ within $X$
and update $m$ to be the new total number of arcs in the digraph $D$. In
terms of this new $m$,
it suffices to show that for any sufficiently  small  real
$\varepsilon>0$, $D$ admits a bipartition $V(D)=V_1\cup V_2$ with
\[
\min\{e(V_{1},V_{2}),e(V_{2},V_{1})\}\geq\Big(\frac{d-1}{2(2d-1)}-\varepsilon\Big)m.
\]

Let $m_{1}=e(X,Y)+e(Y,X)$ and $m_{2}=e(Y)$. Then $m=m_{1}+m_{2}$.
Given a partition $X=X_{1}\cup X_{2}$, let
\[
\theta(X_1,X_2)=\big(e(X_{1},Y)+e(Y,X_{2})\big)-\big(e(X_{2},Y) +e(Y,X_{1})\big)
\]
to be its $gap$.

We choose the partition $X_1, X_2$ such that
$$|\theta(X_1,X_2)| \mbox{ is minimum} .$$
By the symmetry between $X_1$ and $X_2$, we may assume $$\theta(X_1,X_2)\geq 0.$$
Throughout the proof, we write
$\theta:=\theta(X_1,X_2)$, unless it concerns a different
partition of $X$.

\medskip

{\bf Step 2. Extending $X_1,X_2$
to a partition $W_1,W_2 $ of $V(D)$.}

Let $G$ be the underlying graph of $D$, and let $\tau$ be the number of tight components in $G[Y]$.
As $n$ is sufficiently large, we may apply Lemma \ref{main-lemma} and
conclude that $D$ admits a bipartition $V(D)=W_{1}\cup W_{2}$  such that
$X_{i}\subseteq  W_{i}$ for $i=1,2$, and
\begin{align*}
& \min\{e(W_{1},W_{2}),e(W_{2},W_{1})\}
\\ \geq &\frac{1}{2}\min\{e(X_{1},Y)+e(Y,X_{2}),e(X_{2},Y) +e(Y,X_{1})\}+\frac{e(Y)}{4}+\frac{n-\tau}{8}-\varepsilon n
\\ = &\frac{m_{1}-\theta}{4}+\frac{m_{2}}{4}+\frac{n-\tau}{8}-\varepsilon n
\\ =&\frac{m-\theta}{4}+\frac{n-\tau}{8}-\varepsilon n
\\
  \geq&\Big(\frac{d-1}{2(2d-1)}-\varepsilon\Big)m+\frac{1}{4}\Big(\frac{n}{2}+\frac{m}{2d-1}-\theta-\frac{\tau}{2}\Big)
        \quad (\mbox{since $m\ge dn$}).
\end{align*}

Thus, to prove Theorem \ref{main1}, it suffices to show
\begin{equation}\label{mian-ineq-pf}
\frac{n}{2}+\frac{m}{2d-1}\geq\theta+\frac{\tau}{2}.
\end{equation}

\medskip

{\bf Step 3. Bounding $\theta, m$ and $m_2$.}

Since $\tau\leq n$,  \eqref{mian-ineq-pf} holds if
$\theta\leq\frac{m}{2d-1}$. Thus,  we may assume
\begin{equation}\label{bound-theta-m}
\theta>\frac{m}{2d-1}.
\end{equation}
Therefore, since $\theta\leq n'=|Y|$ (by Lemma \ref{hou1}),
\begin{equation}\label{upper-low-bound-m-easy}
dn\leq m<(2d-1)n'.
\end{equation}

Note that we may assume
\begin{equation}\label{bound-m2-easy}
m_{2}>\frac{d+1}{2(2d-1)}m.
\end{equation}
For, otherwise, taking the sum of outdegrees (respectively,
indegrees) of all vertices in $Y$, we have
\[
\min\{e(X,Y),e(Y,X)\}=\min\Big\{\sum_{v\in Y}d^+(v),\sum_{v\in Y}d^-(v)\Big\}-m_2\geq dn'-m_{2}>\frac{d-1}{2(2d-1)}m,
\]
where the last inequality holds since $m<(2d-1)n'$ (by \eqref{upper-low-bound-m-easy}).
Thus, $V(D)=X\cup Y$ gives the desired partition of $D$.

\medskip

{\bf Step 4. Analyzing small tight components.}

We need to analyze   tight components with small order to obtain better bounds on $m$
and $\tau$.
For $i\ge 1$, let $\mathcal{T}_i$ be the collections of tight
components of order $i$ in $G[Y]$, and let $A_i$ be the set of
those vertices of $D$ each of which is contained  in some member of $\mathcal{T}_i$.
Note that by Lemma
\ref{tight-graph-clique}, each $T\in \mathcal{T}_1$ is an isolated
vertex in $G[Y]$ and each $T\in \mathcal{T}_3$ is a triangle in
$G[Y]$.

We bound the number of arcs between $A_i$ and $X$ in both directions
for $i=1, 3$.
For each $x\in A_1$, $\min\{e(x, X), e(X, x)\}\ge d$ since  $\delta^0(D)\ge d$. Thus,
\begin{equation}\label{bound-A1-X}
\min \big\{ e(A_1, X), e(X, A_1)\big \} \ge d|A_1|.
\end{equation}
For each triangle $T\in \mathcal{T}_3$, $D[V(T)]$ has at most 6 arcs
as $D$ has no parallel arcs. Hence,
\begin{align*}
& \min\{e(V(T), X), e(X, V(T))\}
\\= & \min\Big\{\sum_{x\in V(T)}d^+(x), \sum_{x\in V(T)}d^-(x)\Big\}-e(D[V(T)])
\\  \ge & (d-2)|T|.
\end{align*}
Note that the last inequality holds as $\delta^0(D)\ge d$ and $|T|=3$. So
\begin{equation}\label{bound-A3-X}
\min\{e(A_3, X), e(X, A_3)\}\ge \sum_{T\in \mathcal{T}_3}(d-2)|T|=(d-2)|A_3|.
\end{equation}
\medskip

{\bf Step 5. Better bounds on $\tau$ and $m$.}

First, we find a good partition of $D':=D[Y]-(A_1\cup A_3)$. Let $V(D')=U_1\cup U_2$ be a (random) partition of $V(D')$
obtained by  placing each  $v\in V(D')$ into $U_1$ or $U_2$,
independently, with probability $1/2$.
For each arc $e = uv\in A(D')$, let $I_e$ be the indicator random
variable of the event that $u\in U_1$ and $v\in U_2$.
Then, by the linearity of expectation,
\[
\mathbb{E}[e_{D'}(U_1, U_2)]=\sum_{e\in A(D')}\mathbb{E}[I_e]=\sum_{e\in A(D')}\Big(\mathbb{P}[u\in U_1]\mathbb{P}[v\in U_2]\Big)=\frac 1 4e(D').
\]
Similarly,
\[
\mathbb{E}[e_{D'}(U_2, U_1)]=\frac 1 4e(D').
\]
Note  that changing the placement of each  $v\in V(D')$ cannot
affect $ e_{D'}(U_1,  U_2)$ or  $e_{D'}(U_2, U_1)$ by more than
$d_{D'}(v)\le n^{3/4}$; so we may apply Lemma \ref{AH}. Note
\[
L:=\sum_{v\in V(D')}d^2_{D'}(v)\le n^{3/4}\sum_{v\in V(D')}d_{D'}(v)=2n^{3/4}e(D').
\]
Let $z=2n^{3/8}\big(e(D')\big)^{1/2}$. Then by Lemma \ref{AH},
\begin{equation*}\label{py1}
  \mathbb{P}\Big(e_{D'}(U_1, U_2)\leq \mathbb{E}[e_{D'}(U_1, U_2)]-z\Big)\leq \exp\Big(-\frac{z^2}{2L}\Big)= e^{-1},
\end{equation*}
and (by symmetry)
\begin{equation*}\label{py1}
  \mathbb{P}\Big(e_{D'}(U_2, U_1)\leq \mathbb{E}[e_{D'}(U_2, U_1)]-z\Big)\leq  e^{-1}.
\end{equation*}
Thus, with positive probability, there exists a partition $V(D')=U_1\cup U_2$ such that
\begin{equation}\label{bound-w1-w2}
\min\{e_{D'}(U_1,U_2),e_{D'}(U_2,U_1)\}\ge \frac 1 4e(D')-2n^{3/8}\big(e(D')\big)^{1/2}=\frac 1 4e(D')-O(m^{7/8}),
\end{equation}
where the equality holds because of \eqref{m=on=delta}.

\medskip

We now derive better bounds on  $m$ and $\tau$, by working with the
following partition of $D$:  $V_1=U_1\cup A_1\cup A_3$ and $V_2=X\cup U_2$.
Note that $e(A_i, V(D'))=e(V(D'), A_i)=0$ for $i=1,3$. Also note that, for each $T\in
\mathcal{T}_3$, $e(D[V(T)])\le 6$; so
$e(D')\ge m_2-2|A_3|$.
Hence,
\begin{align*}
e(V_1, V_2) & \ge e(A_1, X)+e(A_3,X)+e_{D'}(U_1, U_2)
\\& \ge d|A_1|+(d-2)|A_3|+\frac 1 4e(D')-O(m^{7/8}) \quad (\mbox{by \eqref{bound-A1-X}, \eqref{bound-A3-X} and
\eqref{bound-w1-w2}})
\\ & \ge \frac 1 4m_2+d|A_1|+\left(d-\frac 5 2\right)|A_3|-O(m^{7/8})
     \quad \mbox{(since $e(D')\ge m_2-2|A_3|$)}.
\end{align*}
Similarly,
\begin{align*}
e(V_2, V_1) & \ge e(X, A_1)+e(X, A_3)+e_{D'}(U_2, U_1)
\\ & \ge \frac 1 4m_2+d|A_1|+\left(d-\frac 5 2\right)|A_3|-O(m^{7/8}).
\end{align*}
We may assume that
\begin{equation}\label{assume-w1-bipa-m}
\frac 1 4m_2+d|A_1|+\left(d-\frac 5 2\right)|A_3|\le \frac{d-1}{2(2d-1)}m,
\end{equation}
since, otherwise, $V(D)=V_1\cup V_2$ gives the desired bipartition for
$D$ (as $n$ and, hence, $m$ are sufficiently large).
Thus,
$$|A_1| \le \frac{1}{d}\left(\frac{d-1}{2(2d-1)}m-\frac 1 4m_2-\left(d-\frac 5 2\right)|A_3|\right)
 < \frac {3d-5}{8d(2d-1)}m-\left(1-\frac{5}{2d}\right)|A_3|,$$
since $m_2>\frac {d+1}{2(2d-1)}m$ (by \eqref{bound-m2-easy}). Therefore,
\begin{align}\label{now-bound-tau}
\tau & \le |A_1|+\frac{|A_3|}{3}+\frac{n'-|A_1|-|A_3|}{5} \notag
\\ & =\frac{n'}{5}+\frac{4}{5}|A_1|+\frac{2}{15}|A_3|         \notag
\\ & \le \frac{n'}{5}+\frac{4}{5}\left(\frac {3d-5}{8d(2d-1)}m-\left(1-\frac{5}{2d}\right)|A_3|\right)+\frac{2}{15}|A_3|   \notag
\\ & =\frac{1}{5}n'+\frac{3d-5}{10d(2d-1)}m-\Big(\frac{2}{3}-\frac{2}{d}\Big)|A_3|   \notag
\\ & < \frac{1}{5}n'+\frac{3d-5}{10d(2d-1)}m \quad \mbox{(as $d\ge 4$)}.
\end{align}
Then
\begin{align*}
& \frac{n}{2}+\frac{m}{2d-1}-\theta-\frac{\tau}{2}
\\ \ge &
         \frac{n'}{2}+\frac{m}{2d-1}-\theta-\frac{1}{2}\left(\frac{1}{5}n'+\frac{3d-5}{10d(2d-1)}m\right)
         \quad \mbox{ (by \eqref{now-bound-tau})}
\\ =&  \frac{2}{5}n'+\frac{17d+5}{20d(2d-1)}m-\theta
\\ \ge &  \frac{17d+5}{20d(2d-1)}m-\frac{3}{5}n' \quad \mbox{ (since
         $\theta\le n'$)}.
\end{align*}
Thus, if $m\ge \frac{12d(2d-1)}{17d+5}n'$ then  \eqref{mian-ineq-pf} holds and, hence,  $W_1,W_2$ give the desired
partition for $D$.
So we may assume that
\begin{equation}\label{now-bound-m}
m< \frac{12(2d-1)}{17d+5}dn'.
\end{equation}

\medskip

{\bf Step 6. Bounding $m_1$.}

To establish an upper bound on  $m_1$,
we  sum up the  outdegrees
(respectively, indegrees) of all vertices in $Y$. Thus, by $\delta^0(D)\ge d$,
$\min\{e(X,Y),e(Y,X)\}\geq dn'-m_{2}$; so by \eqref{now-bound-m},
\[
\min\{e(X,Y),e(Y,X)\}\geq \frac{17d+5}{12(2d-1)}m-m_2.
\]

If $m_2\le \frac{11d+11}{12(2d-1)}m$ then
\[
\min\{e(X,Y),e(Y,X)\}\geq \frac{d-1}{2(2d-1)}m,
\]
and $X,Y$ form the desired partition for $D$.

So we may assume
$m_2> \frac{11d+11}{12(2d-1)}m$. Then by \eqref{now-bound-m}, we have
\begin{equation}\label{now-bound-m1}
m_1=m-m_2< \frac{13d-23}{12(2d-1)}m< \frac{13d-23}{17d+5}dn'.
\end{equation}

\medskip

{\bf Step 7. Bounding $\alpha$,  the number of huge vertices}.

  For a vertex $v\in X$,  let
$s^- (v) = d^-(v)-d^+ (v)$,
$s^+ (v) = d^+ (v) -d^- (v)$,
and $s(v) = \max\{s^+ (v),s^-(v)\}$.
As in \cite{Lee2014}, we call a vertex $v\in X$ $huge$ if $s(v)\geq\theta$.
For the partition  $X = X_1\cup X_2$ (previously chosen to minimize
$\theta$) and  a vertex $v\in X$, $v$  is \emph{forward} if either $v\in X_1$ and $s^+ (v) > 0$,
or $v \in X_2$ and $s^- (v) > 0$. Since $\theta>0$ (by \eqref{bound-theta-m}), $X$ has at
least one forward vertex.

We claim that all forward vertices in $X$ are huge.
For, let $v$ be a forward vertex in $X$, and let the partition
$X_1'\cup X_2'$ of $X$ be obtained from $X_1\cup X_2$ by switching
the side of $v$. Then
\begin{align*}
\theta(X'_1, X'_2) & =\big(e(X'_1,Y)+e(Y,X'_2)\big)-\big(e(X'_2,Y) +e(Y,X'_1)\big)
\\& = \big(e(X_1,Y)+e(Y,X_2)-s(v)\big)-\big(e(X_2,Y) +e(Y,X_1)+s(v)\big)
\\ & = \theta-2s(v).
\end{align*}
Since $X=X_1\cup X_2$ is   a partition of $X$  minimizing $\theta=|\theta(X_1,X_2)|$,
$|\theta(X'_1, X'_2)|\ge \theta$. Therefore,
$\theta-2s(v)\leq-\theta$; so  $s(v)\geq \theta$ and, hence, $v$ is huge.

Let $X'$ be the set of huge vertices in $X$ and $X''=X- X'$. Let
$\alpha:=|X'|$ and $\rho:=\sum_{v\in X''}{d(v)}$. Then $\alpha \ge 1$,
and
\begin{equation}\label{bound-m1-easy}
m_{1}=\sum_{v\in X}{d(v)}=\sum_{v\in X'}{d(v)}+\sum_{v\in X''}{d(v)}\geq\alpha\theta+\rho.
\end{equation}
Since $\theta>m/(2d-1)$ (by \eqref{bound-theta-m}), it follows from
\eqref{now-bound-m1} and \eqref{bound-m1-easy} that
\[
\frac{\alpha}{2d-1}m\le \frac{13d-23}{12(2d-1)}m,
\]
which implies
\begin{equation}\label{1-f(alpha)}
1-\frac{\alpha}{2(2d-1)}>0.
\end{equation}

\medskip

To derive a better upper bound on $\alpha$,
we sum up the  degrees of all vertices in $Y$. So
$$2dn'\leq\sum_{v\in Y}{d(v)}=e(X,Y)+e(Y,X)+2e(Y)=m_{1}+2m_{2},$$
which, together with \eqref{bound-m1-easy}, gives
\begin{equation}\label{another-bound-m}
m=m_{1}+m_{2}\geq\frac{2dn'+m_{1}}{2}\geq\frac{2dn'+\alpha\theta+\rho}{2}.
\end{equation}
By \eqref{now-bound-tau} and \eqref{now-bound-m}, we have
\[
\tau\le \frac{1}{5}n'+\frac{3d-5}{10d(2d-1)}
\frac{12d(2d-1)}{17d+5}n'=\frac{7d-5}{17d+5}n';
\]
so
 \begin{align*}
& \frac{n}{2}+\frac{m}{2d-1}-\theta-\frac{\tau}{2}  \notag
\\ \ge & \frac{n}{2}+\frac{m}{2d-1}-\theta-\frac{7d-5}{2(17d+5)}n' \notag
\\ \ge &  \frac{n'}{2}+\frac{1}{2(2d-1)}\Big(2dn'+\alpha\theta+\rho\Big)-\theta-\frac{7d-5}{2(17d+5)}n'      \quad (\mbox{by \eqref{another-bound-m}})  \notag
\\ \ge & \left(\frac{1}{2}+\frac{2d}{2(2d-1)}-\frac{7d-5}{2(17d+5)}\right)n'-\left(1-\frac{\alpha}{2(2d-1)}\right)\theta       \notag
\\ \ge &
         \Big(\frac{1}{2}+\frac{2d}{2(2d-1)}-\frac{7d-5}{2(17d+5)}\Big)n'-\Big(1-\frac{\alpha}{2(2d-1)}\Big)n'      \quad (\mbox{by \eqref{1-f(alpha)} and $\theta\le n'$}) \notag
\\ = & \Big(\frac{\alpha+1}{2(2d-1)}-\frac{7d-5}{2(17d+5)}\Big)n'. \notag
\end{align*}
Hence, if $\alpha\ge \frac{2d(7d-17)}{17d+5}$ then
$$\frac{n}{2}+\frac{m}{2d-1}-\theta-\frac{\tau}{2}\ge \Big(\frac{\alpha+1}{2(2d-1)}-\frac{7d-5}{2(17d+5)}\Big)n'>0;$$
so \eqref{mian-ineq-pf} holds,  and
$W_1,W_2$ gives the desired partition. We may thus assume that
\begin{equation}\label{small-alpha}
\alpha< \frac{2d(7d-17)}{17d+5}<d-1.
\end{equation}

\medskip
{\bf Step 8. Better bounds on $\tau$.}

We bound $\tau$  in terms of $n'$ and $\rho$. Let $d_{i}$ denote the number of vertices with degree $i$ in $G[Y]$.  Then
\[
\sum_{i\ge0}d_i= n'.
\]
For $0\le i\le d-\alpha-1$, each vertex $v\in Y$ with degree $i$ in
$G[Y]$ has degree at most $2i$ in $D[Y]$, which, together with the fact
that $\delta^0(D)\ge d$, shows  $e(v, X)+e(X, v)\ge 2d-2i$.  Thus, the
number of arcs of $D$ incident with $v$ and counted in $\rho$ is at
least $2d-2i-2\alpha$.
Consequently,
\begin{equation*}\label{tight-bound-rho}
(2d-2\alpha)d_0+(2d-2\alpha-2)d_1+\cdots+2d_{d-\alpha-1}\leq\rho,
\end{equation*}
i.e.,
\begin{equation*}\label{tight-bound-rho}
(d-\alpha)d_0+(d-\alpha-1)d_1+\cdots+d_{d-\alpha-1}\leq\rho/2.
\end{equation*}
Therefore, by Lemmas \ref{Lee1} and \ref{hou2}, we have
\begin{equation}\label{final-bound-tau}
\tau\leq\frac{n'+\rho/2}{d-\alpha+1}.
\end{equation}

\medskip

{\bf Step 9. Completing the proof by considering $\alpha =1$ and
  $\alpha\ge 2$.}

For convenience, let $m_1=\beta dn'$. Then, by \eqref{now-bound-m1}, we have
\begin{equation}\label{final-bound-beta}
\beta\le \frac{13d-23}{17d+5}.
\end{equation}
A simple calculation shows that
\begin{equation}\label{ineq-beta}
(1-\beta)d-\frac{3}{2}> 0.
\end{equation}

If $\alpha=1$ then
\begin{align*}
& \frac{n}{2}+\frac{m}{2d-1}-\theta-\frac{\tau}{2}  \\
{\ge} &
        \frac{n'}{2}+\frac{1}{2(2d-1)}\Big(2dn'+\theta+\rho\Big)-\theta
        -\frac{n'+\rho/2}{2d}  \quad (\mbox{ by
        \eqref{another-bound-m} and \eqref{final-bound-tau}})\\
= & \Big(\frac{1}{2}+\frac{2d}{2(2d-1)}-\frac{1}{2d}\Big)n'-\Big(1-\frac{1}{2(2d-1)}\Big)\theta +\Big(\frac{1}{2(2d-1)}-\frac{1}{4d}\Big)\rho \\
\ge &
      \Big(\frac{1}{2}+\frac{2d}{2(2d-1)}-\frac{1}{2d}\Big)n'-\Big(1-\frac{1}{2(2d-1)}\Big)n'
      \quad (\mbox{since $\theta\le n'$ by Lemma \ref{hou1}})                \\
 = & \Big(\frac{1}{2d-1}-\frac{1}{2d}\Big)n'\\
 > & 0.
\end{align*}
Hence,  \eqref{mian-ineq-pf} holds, and
$W_1,W_2$ gives the desired partition.

\medskip

So we may assume $\alpha \ge 2$ and, hence,
\begin{equation}\label{alpha-ineq1}
\frac{1}{4(d-\alpha+1)}>\frac{1}{2(2d-1)}.
\end{equation}
Note that
\begin{align*}
& \frac{n}{2}+\frac{m}{2d-1}-\theta-\frac{\tau}{2}  \notag
\\
\ge &  \frac{n'}{2}+\frac{1}{2(2d-1)}\Big(2dn'+\alpha\theta+\rho\Big)-\theta -\frac{n'+\rho/2}{2(d-\alpha+1)}    \quad (\mbox{by
        \eqref{another-bound-m} and \eqref{final-bound-tau}})  \notag
\\ = & \Big(\frac{1}{2}+\frac{2d}{2(2d-1)}-\frac{1}{2(d-\alpha+1)}\Big)n'-\Big(1-\frac{\alpha}{2(2d-1)}\Big)\theta       \notag
\\     & -\Big(\frac{1}{4(d-\alpha+1)}-\frac{1}{2(2d-1)}\Big)\rho
\\
>&
      \Big(\frac{4d-1}{2(2d-1)}-\frac{1}{2(d-\alpha+1)}\Big)n'-\Big(1-\frac{\alpha}{2(2d-1)}\Big)\theta
\\     & -\Big(\frac{1}{4(d-\alpha+1)}-\frac{1}{2(2d-1)}\Big)\Big(\beta dn'-\alpha\theta\Big)
 \quad (\mbox{by \eqref{bound-m1-easy}, \eqref{alpha-ineq1}
      and since $m_1=\beta dn'$})   \notag
\\ =&    \Big(\frac{(4+\beta)d-1}{2(2d-1)}-\frac{\beta d+2}{4(d-\alpha+1)}\Big)n'-\Big(1-\frac{\alpha}{4(d-\alpha+1)}\Big)\theta.                   \notag
\end{align*}

First, suppose $1-\frac{\alpha}{4(d-\alpha+1)}\ge 0$. Then
\begin{align*}
& \frac{n}{2}+\frac{m}{2d-1}-\theta-\frac{\tau}{2}  \notag
\\ \ge & \Big(\frac{(4+\beta)d-1}{2(2d-1)}-\frac{\beta
         d+2}{4(d-\alpha+1)}\Big)n'-\Big(1-\frac{\alpha}{4(d-\alpha+1)}\Big)n'
         \quad (\mbox{since $\theta\le n'$}) \notag
\\  = & \frac{(2d-1)(\alpha-1)-(\beta d+1)(2\alpha-3)}{4(2d-1)(d-\alpha+1)}n'
\\  > & \frac{(2\alpha-3)\big((1-\beta)d-3/2\big)}{4(2d-1)(d-\alpha+1)}n'
\\ >& 0 \quad (\mbox{since $\alpha \ge 2$
                                  and by \eqref{ineq-beta}}).
\end{align*}
So \eqref{mian-ineq-pf} holds, and
$W_1,W_2$ gives the desired partition.

Thus, we may assume
$1-\frac{\alpha}{4(d-\alpha+1)}<0$; so $\alpha>4d/5$.
Note that
\begin{align*}
& \frac{n}{2}+\frac{m}{2d-1}-\theta-\frac{\tau}{2}  \notag
\\ > & \Big(\frac{(4+\beta)d-1}{2(2d-1)}-\frac{\beta d+2}{4(d-\alpha+1)}\Big)n'-\Big(1-\frac{\alpha}{4(d-\alpha+1)}\Big)\frac{m}{2d-1}
\quad (\mbox{by \eqref{bound-theta-m}})
\\ \ge &  \Big(\frac{(4+\beta)d-1}{2(2d-1)}-\frac{\beta d+2}{4(d-\alpha+1)}\Big)n'-\Big(1-\frac{\alpha}{4(d-\alpha+1)}\Big)\frac{d}{2d-1}n
\quad (\mbox{by the fact $m\ge dn$})
\\ = &  \Big(\frac{(2+\beta)d-1}{2(2d-1)}-\frac{\beta
       d+2}{4(d-\alpha+1)}+\frac{\alpha
       d}{4(d-\alpha+1)(2d-1)}\Big)n'-O(n^{1/4}) \quad (\mbox{as
       $n'=n-O(n^{1/4})$}) \notag
\\ = &  \Big(\frac{4d^2-3\alpha d-2\alpha\beta d-2d+2\alpha+3\beta d}{4(d-\alpha+1)(2d-1)}\Big)n'-O(n^{1/4})
\\ > &  \Big(\frac{4d-3\alpha -2\alpha\beta
         -2/5}{4(d-\alpha+1)(2d-1)}\Big)dn'-O(n^{1/4}) \quad
         (\mbox{since $\alpha >4d/5$}) \notag
\\ > &  \Big(\frac{\big(4-3\cdot\frac{14}{17} -2\cdot\frac{14}{17}\cdot\frac{13}{17}\big)d-2/5}{4(d-\alpha+1)(2d-1)}\Big)dn'-O(n^{1/4})
\quad (\mbox{$\alpha<\frac{14d}{17}$ by \eqref{small-alpha}, $\beta<\frac{13}{17}$ by \eqref{final-bound-beta}) }
\\ > &  \Big(\frac{d/4-2/5}{4(d-\alpha+1)(2d-1)}\Big)dn'-O(n^{1/4})
\\ >& 0 \quad (\mbox{since $d\ge 4$ and $n'=n-O(n^{1/4})$ is large}).
\end{align*}
Again, \eqref{mian-ineq-pf} holds, and
$W_1,W_2$ gives the desired partition.

\section{Concluding remarks}

Using copies of $K_{2d-1}$ and one copy of $K_{2d+1}$, Lee, Loh, and Sudakov \cite{Lee2014} constructed digraphs with  minimum
outdegree $d$ and  minimum indegree $d-1$ to show that the main term $\frac{d-1}{2(2d-1)}$ in Conjecture \ref{conjecture} is best
possible. One could even ask whether or not Conjecture \ref{conjecture} holds
without the $o(1)$ term, and, in particular, for $d=2,3$ (as for these
cases
Conjecture~\ref{conjecture} is known to be true). More precisely, is it true that if $D$ is a
digraph with $m$ arcs and the minimum outdegree at
least 2 (respectively, 3) then $D$ admits a bipartition $V(D)=V_1\cup V_2$ such that
$\min\{e(V_1,V_2),e(V_2,V_1)\}\ge m/6$ (respectively, $\min\{e(V_1,V_2),e(V_2,V_1)\}\ge m/5$)?

It is not clear to us whether or not the minimum outdegree
condition alone is sufficient for Conjecture~\ref{conjecture} with
large $d$. We have demonstrated that the assertion of Conjecture \ref{conjecture}
holds if we impose  the additional (natural) condition that the minimum indegree of
$D$ be  also at least $d$. (Though we do not know whether the main term in Theorem \ref{main1} is tight.)
Below, we comment on the five places in our proof where we
use this indegree assumption.

The first place is where we prove
\eqref{bound-m2-easy} for a lower bound on $m_2$ (the number of arcs in
$D[Y]$). The minimum
outdegree condition alone is not sufficient for this purpose. However,
we checked several extreme situations (e.g., when $e(Y)=0$), we can
always find a partition satisfying the conjectured bound. But we do not know how to deal with the general case.

The second place is in the proof of \eqref{bound-A1-X} and
\eqref{bound-A3-X}. We need to bound the number of arcs of $D$  between
$X$ and $A_i$ in both directions  for $i=1,3$, in order to bound
$\tau$ (the number of tight components in $G[Y]$) in
\eqref{now-bound-tau} and  $m$ (the number of arcs in
$D$) in \eqref{now-bound-m}.

The third place is where we
obtain $\min\{e(X,Y),e(Y,X)\}\ge dn'-m_2$  to bound $m_1$ in \eqref{now-bound-m1}, and the fourth place is where
we prove the lower bound on $m$ in \eqref{another-bound-m}.
The final place we use this minimum indegree condition is where we
give an upper bound on $\tau$ in \eqref{final-bound-tau}.



\end{document}